\documentclass[11pt]{article}

\usepackage{graphicx,psfrag}
\usepackage{psfrag,eepic,epsfig}
\usepackage{sectsty}
\usepackage{setspace}
\usepackage{color}

\usepackage{setspace}
\usepackage{amsmath, epsfig, cite, ifpdf}
\usepackage{amssymb}
\usepackage{amsfonts}
\usepackage{latexsym}
\usepackage{graphicx,psfrag,eepic,rawfonts}
\usepackage{enumerate}
\usepackage{indentfirst}
\input{epsf}

\newtheorem{thm}{Theorem}[section]

\newtheorem{defi}[thm]{Definition}

\newtheorem{lem}[thm]{Lemma}

\newtheorem{problem}[thm]{Problem}

\newcommand{\qed}{{\hfill\rule{4pt}{7pt}}}
\def\pf{\noindent {\it Proof.} }

\numberwithin{equation}{section}

\makeatletter \@addtoreset{equation}{section} \makeatother

\title {\bf $3$-Regular mixed graphs with \\optimum Hermitian
energy\footnote{Supported by NSFC No.11371205, the 973 program of China No.2013CB834204, and PCSIRT. }}

\author{
{\small Xiaolin Chen, Xueliang Li, Yingying Zhang}\\
{\small Center for Combinatorics and LPMC-TJKLC}\\
{\small Nankai University, Tianjin 300071, P.R. China}\\
{\small E-mail: chxlnk@163.com; lxl@nankai.edu.cn; zyydlwyx@163.com}
   }
\date{}
\begin{document}

\maketitle

\begin{abstract}
Let $G$ be a simple undirected graph, and $G^\phi$ be a mixed graph
of $G$ with the generalized orientation $\phi$ and
Hermitian-adjacency matrix $H(G^\phi)$. Then $G$ is called the
underlying graph of $G^\phi$. The Hermitian energy of the mixed
graph $G^\phi$, denoted by $\mathcal{E}_H(G^\phi)$, is defined as
the sum of all the singular values of $H(G^\phi)$. A $k$-regular
mixed graph on $n$ vertices having Hermitian energy $n\sqrt{k}$ is
called a $k$-regular optimum Hermitian energy mixed graph. In this
paper, we first focus on the problem proposed by Liu and Li [J. Liu,
X. Li, Hermitian-adjacency matrices and Hermitian energies of mixed
graphs, Linear Algebra Appl. 466(2015), 182--207] of determining all
the $3$-regular connected optimum Hermitian energy mixed graphs. We
then prove that optimum Hermitian energy oriented graphs with
underlying graph hypercube are unique (up to switching equivalence).


\end{abstract}

\section{Introduction}

Let $G$ be a simple undirected graph with vertex set
$V(G)$ and edge set $E(G)$. A generalized orientation $\phi$ of $G$ is to give each edge of $S$ an
orientation according to $\phi$, where $S\subseteq{E(G)}$. Then $G^{\phi}$ is called a mixed graph
of $G$ with the generalized orientation $\phi$. If $S={E(G)}$, $\phi$ is an orientation of $G$ and the mixed graph $G^{\phi}$ is an oriented
graph. If $S=\emptyset$, then $G^{\phi}$ is an undirected graph. Thus we find that mixed graphs incorporate
both undirected graphs and oriented graphs as extreme cases. In a mixed graph $G^{\phi}=(V(G^{\phi}),E(G^{\phi}))$,
if one element $(u,v)$ in $E(G^{\phi})$ is an edge (resp. arc), we denote it by $u\leftrightarrow{v}$
(resp. $u\rightarrow{v}$). The graph $G$ is called the underlying graph of $G^{\phi}$.
A mixed graph is called regular if its underlying graph is a regular graph. Similarly, in terms of defining
order, size, degree and so on, we focus only on its underlying graph. For undefined terminology and notations,
we refer the reader to \cite{6,9}.

The Hermitian-adjacency matrix $H(G^{\phi})$ of $G^{\phi}$ with vertex set $V(G^{\phi})=\{1,2,\ldots,n\}$ is
a square matrix of order n, whose entry $h_{kl}$ is defined as
\begin{equation*}
h_{kl}=
\begin{cases}
h_{lk}=1, & if \ \ k\leftrightarrow{l},\\
-h_{lk}=i, & if \ \ k\rightarrow{l},\\
\ \ 0, & otherwise,
\end{cases}
\end{equation*}
where $i$ is the imaginary number unit. The spectrum
$Sp_{H}(G^{\phi})$ of $G^{\phi}$ is defined as the spectrum of
$H(G^{\phi})$. Since $H(G^{\phi})$ is a Hermitian matrix, i.e.,
$H(G^{\phi})=[H(G^{\phi})]^{*}:=\overline{[H(G^{\phi})]}^{T}$, the
eigenvalues $\{\lambda_{1},\lambda_{2},\ldots,\lambda_{n}\}$ of
$H(G^{\phi})$ are all real. In \cite{1}, Liu and Li introduced the
Hermitian energy of the mixed graph $G^{\phi}$, denoted by
$\mathcal{E}_H(G^{\phi})$, which is defined as the sum of the
singular values of $H(G^{\phi})$. Since the singular values of
$H(G^{\phi})$ are the absolute values of its eigenvalues, we have
$$\mathcal{E}_H(G^{\phi})=\sum_{j=1}^{n}|\lambda_{j}|.$$

For an oriented graph $G^{\phi}$, Adiga et al. \cite{ABS} introduced
the concept of skew adjacency matrix of $G^{\phi}$, denoted by
$S(G^{\phi})$, which is defined as $S(G^{\phi})=-iH(G^{\phi})$.
Then, the eigenvalues of $S(G^{\phi})$ are
$\{-i\lambda_{1},-i\lambda_{2},\ldots,-i\lambda_{n}\}$. The skew
energy of oriented graph $G^{\phi}$ is defined by Adiga et al. in
\cite{ABS} as
$\mathcal{E}_S(G^{\phi})=\sum_{j=1}^{n}|-i\lambda_{j}|$. Thus,
$\mathcal{E}_H(G^{\phi})=\mathcal{E}_S(G^{\phi})$, i.e., the
Hermitian energy of an oriented graph is equal to its skew energy.
For more details about skew energy, we refer the survey \cite{LL} to
the reader.

Hermitian energy can be viewed as a generalization of the graph
energy. The concept of the energy of simple undirected graphs was
introduced by Gutman in \cite{7}, which is related to the total
$\pi$-electron energy of the molecule represented by that graph.
Since then, the graph energy has been extensively studied. For more
details, we refer \cite{12} to the reader.

In \cite{1}, Liu and Li gave a sharp upper bound of the Hermitian
energy in terms of its order $n$ and the maximum degree $\Delta$,
i.e. $$\mathcal{E}_H(G^{\phi})\leq n\sqrt{\Delta}.$$ Furthermore,
they showed that the equality holds if and only if
$H^{2}(G^{\phi})=\Delta{I_{n}}$, which implies that $G^{\phi}$ is
$\Delta$-regular. For convenience, in this paper a mixed graph on
$n$ vertices with maximum degree $\Delta$ which satisfies
$\mathcal{E}_H(G^{\phi})=n\sqrt{\Delta}$ is called an optimum
Hermitian energy mixed graph. Let $I_{n}$ be the identity matrix of
order $n$. For simplicity, we always write $I$ when its order is
clear from the context. It is important to determine a family of
$k$-regular mixed graphs with optimum Hermitian energy for any
positive integer $k$. In \cite{1}, Liu and Li gave $Q_{k}$ a
suitable generalized orientation such that it has optimum Hermitian
energy. Besides, they proposed the following problem:
\begin{problem}\label{PRO}
Determine all the $k$-regular mixed graphs $G^{\phi}$ on $n$
vertices with $\mathcal{E}_H(G^{\phi})=n\sqrt{k}$ for each $k$,
$3\leq{k}\leq{n}$.
\end{problem}

Liu and Li \cite{1} showed that a $1$-regular connected mixed graph
on $n$ vertices has optimum Hermitian energy if and only if it is an
edge or arc. At the same time, they also proved that a $2$-regular
connected mixed graph on $n$ vertices has optimum Hermitian energy
if and only if it is one of the three types of mixed $4$-cycles. If
$G_{1}^{\phi}$ and $G_{2}^{\phi}$ are two $k$-regular mixed graphs
with optimum Hermitian energy, then so is their disjoint union.
Thus, we only consider $k$-regular connected mixed graphs.

In this paper, we firstly characterize all $3$-regular connected
optimum Hermitian energy mixed graphs. Thus we solve Problem
\ref{PRO} for $k=3$. Afterwards, we prove that optimum Hermitian
energy oriented graphs with underlying graph hypercube are unique
(up to switching equivalence).

\section{Preliminaries}

In this section, we give some notations and known results. Besides, we also introduce the
definition of switching equivalence.

Let $G=G(V,E)$ be a graph with vertex set $V$ and edge set $E$. For any $v\in{V}$,
we denote the neighborhood of $v$ by $N_{G}(v)$ in $G$. Let $G[S]$ denote the subgraph
of $G$ induced by $S$, where $S\subseteq{V}$. In addition, we give $G$ a generalized
orientation $\phi$. Then we get a mixed graph denoted by $G^\phi=(V(G^\phi),E(G^\phi))$
and the Hermitian-adjacency matrix of $G^\phi$ by $H(G^\phi)$.

In \cite{1}, Liu and Li gave a sharp upper bound for the Hermitian energy of a mixed graph
and a necessary and sufficient condition to attain the upper bound.

\begin{lem}{(12, a part of Theorem 3.2).}\label{L1}
Let $G^{\phi}$ be a mixed graph on $n$ vertices with maximum degree $\Delta$. Then
$\mathcal{E}_H(G^{\phi})\leqslant{n\sqrt{\Delta}}$.
\end{lem}

\begin{lem}{(12, a part of Corollary 3.3).}\label{L2}
Let $H$ be the Hermitian-adjacency matrix of a mixed graph $G^{\phi}$ on $n$ vertices.
Then $\mathcal{E}_H(G^{\phi})={n\sqrt{\Delta}}$ if and only if $H^{2}=\Delta{I_{n}}$ i.e.
the inner products $H(u,:)\cdot{H(v,:)}=0,\  H(:,u)\cdot{H(:,v)}=0$ for different vertices
$u$ and $v$ of $G^{\phi}$, where $H(u,:)$ and $H(:,u)$ represent row vector and column
vector corresponding to vertex $u$ in $H(G^{\phi})$, respectively.
\end{lem}

\begin{figure}[h,t,b,p]
\begin{center}
\scalebox{1}[1]{\includegraphics{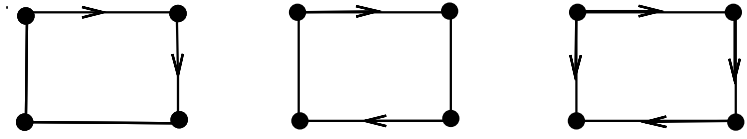}}
\end{center}
\caption{Three types of mixed $4$-cycles.}\label{Fig.2.1.}
\end{figure}

Moreover, Liu and Li \cite{1} gave a characterization of the $k$-regular connected optimum
Hermitian energy mixed graphs.
\begin{lem}{(12, a part of Lemma 3.5).}\label{L3}
Let $G^{\phi}$ be a $k$-regular connected mixed graph with order $n$ $(n\geq{3})$, then
$\mathcal{E}_H(G^{\phi})={n\sqrt{k}}$ if and only if for any pair of vertices $u$ and $v$
with distance not more than two in $G$ such that $N(u)\cap{N(v)}\neq{\emptyset}$, there
are edge-disjoint mixed $4$-cycles $uxvy$ of the following three types; see Fig.\ref{Fig.2.1.}.
\end{lem}

By Lemma \ref{L2}, if $G^{\phi}$ is a connected mixed graph on $n$ vertices with optimum
Hermitian energy $n\sqrt{\Delta}$, then $G^{\phi}$ is $\Delta{}$-regular. Moreover, since
any two distinct rows of $H$ are orthogonal, we deduce the following lemma.

\begin{lem}\label{L4}
Let $H$ be the Hermitian-adjacency matrix of a $k$-regular mixed graph $G^{\phi}$ on $n$
vertices. If $H^{2}=k{I_{n}}$, then $|N(u)\cap{N(v)}|$ is even for any pair of vertices $u$
and $v$ with distance no more than two in $G$.
\end{lem}

Next we introduce the definition of switching equivalence. Let
$G^{\phi}$ be a mixed graph with vertex set $V$. The switching
function of $G^{\phi}$ is a function
$\theta:V\rightarrow{\mathrm{T}}$, where $\mathrm{T}=\{1,-1\}$. The
switching matrix of $G^{\phi}$ is a diagonal matrix
$D(\theta):=diag(\theta{(v_{k})}:v_{k}\in{V})$, where $\theta$ is a
switching function. Let $G^{\phi_{1}},G^{\phi_{2}}$ and
$G^{\phi_{3}}$ be three mixed graphs with the same underlying graph
$G$ and vertex set $V$. If there exists a switching matrix
$D(\theta)$ such that
$H(G^{\phi_{2}})=D(\theta)^{-1}H(G^{\phi_{1}})D(\theta)$, then we
say $G^{\phi_{1}}$ and $G^{\phi_{2}}$ are switching equivalent,
denoted by $G^{\phi_{1}}\sim{G^{\phi_{2}}}$. If two mixed graphs
$G^{\phi_{1}}$ and $G^{\phi_{2}}$ are switching equivalent, then
$Sp_{H}(G^{\phi_{1}})=Sp_{H}(G^{\phi_{2}})$ i.e.
$\mathcal{E}_H(G^{\phi_{1}})=\mathcal{E}_H(G^{\phi_{2}})$. Besides,
the number of arcs (or undirected edges) in $G^{\phi_{1}}$ is equal
to that in $G^{\phi_{2}}$. Moreover, if
$G^{\phi_{1}}\sim{G^{\phi_{2}}}$ and
$G^{\phi_{2}}\sim{G^{\phi_{3}}}$, then
$G^{\phi_{1}}\sim{G^{\phi_{3}}}$.

Note that Liu and Li \cite{1} also introduced the definition of switching equivalence
between mixed graphs. However, $\mathrm{T}$ is $\{1,i,-i\}$ in their definition. Besides,
our definition coincides with the definition of switching equivalence between oriented
graphs which is given in \cite{5} when the mixed graphs are oriented graphs.

Let $G^{\phi}$ be a $k$-regular optimum Hermitian energy mixed
graph. If $G^{\phi}$ is an oriented graph, then the Hermitian energy
of $G^{\phi}$ is equal to its skew energy. In \cite{2}, Gong and Xu
characterized the $3$-regular optimum Hermitian energy oriented
graphs. Moreover, Chen et al. \cite{CLL} and Gong et al. \cite{GXZ}
independently characterized the $4$-regular optimum Hermitian energy
oriented graphs. The following lemma is the result about the
characterization of $3$-regular optimum Hermitian energy oriented
graphs in \cite{2}.

\begin{figure}[h,t,b,p]
\begin{center}
\scalebox{0.8}[0.8]{\includegraphics{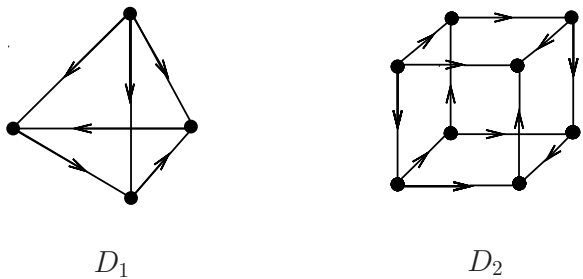}}
\end{center}
\caption{$3$-regular optimum Hermitian energy oriented graphs.}\label{Fig.2.2.}
\end{figure}

\begin{lem}\cite{2}\label{L5}
Let $G^{\phi}$ be a $3$-regular optimum Hermitian energy oriented
graph. Then $G^{\phi}$ (up to isomorphism) is either $D_1$ or $D_2$
drawn in Fig.\ref{Fig.2.2.}.
\end{lem}

\section{The $3$-regular optimum Hermitian energy mixed graphs}

In this section, we characterize all $3$-regular connected optimum
Hermitian energy mixed graphs (up to switching equivalence).

Let $G^{\phi}$ be a $3$-regular optimum Hermitian energy mixed
graph. By Lemma \ref{L4}, we get that $G^{\phi}$ satisfies that
$|N(u)\cap{N(v)}|$ is even for any two distinct vertices $u$ and $v$
of $G^{\phi}$. Moreover based on the proof of Theorem 3.5 in
\cite{2}, we deduce that the underlying graph of $G^{\phi}$ with
$\mathcal{E}_H(G^{\phi})={n\sqrt{3}}$ is either the complete graph
$K_{4}$ or the hypercube $Q_{3}$. Hence, we just need to consider
the $3$-regular optimum Hermitian energy mixed graphs with
underlying graph $K_{4}$ or $Q_{3}$.

Firstly, we consider the case that the underlying graph is the
complete graph $K_4$.

\begin{figure}[h,t,b,p]
\begin{center}
\scalebox{0.8}[0.8]{\includegraphics{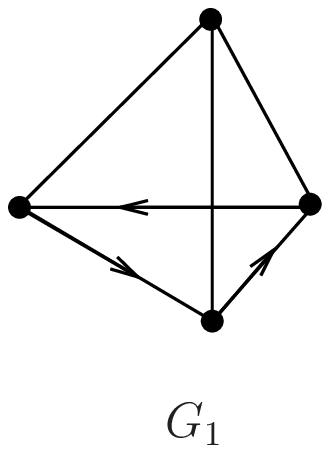}}
\end{center}
\caption{Optimum Hermitian energy mixed graph with underlying graph $K_4$.}\label{Fig.3.3.}
\end{figure}

\begin{thm}\label{T5}
Let $G^{\phi}$ be a $3$-regular optimum Hermitian energy mixed
graph. If the underlying graph $G$ is $K_{4}$, then $G^{\phi}$ is
either $D_{1}$ drawn in Fig.\ref{Fig.2.2.} or $G_{1}$ drawn in
Fig.\ref{Fig.3.3.}.
\end{thm}

\pf We divide our discussion into four cases:

\textbf{Case 1.} $G^{\phi}$ is an oriented graph.

From Lemma \ref{L5}, we obtain
that $G^{\phi}$ is $D_{1}$ drawn in Fig.\ref{Fig.2.2.}.

\textbf{Case 2.} $G^{\phi}$ is not an oriented graph and no vertex has two incident edges.
Then there exists a vertex, say $u_1$, which has one incident edge $u_1\leftrightarrow{u_2}$.

Subcase 2.1. $u_1\rightarrow{u_3}$ and $u_1\rightarrow{u_4}$.

By Lemma \ref{L2}, we have $H(u_1,:)\cdot{H(u_2,:)}=0$. Then $\overline{h}_{11}h_{21}+\overline{h}_{12}h_{22}+\overline{h}_{13}h_{23}+\overline{h}_{14}h_{24}=-ih_{23}-ih_{24}=0$.
Hence $h_{23}=i, h_{24}=-i$ or $h_{23}=-i, h_{24}=i$, that is, $u_{2}\rightarrow{u_{3}},u_{2}\leftarrow{u_{4}}$ or $u_{2}\leftarrow{u_{3}},u_{2}\rightarrow{u_{4}}$. Without loss of generality, assume that $u_{2}\rightarrow{u_{3}},u_{2}\leftarrow{u_{4}}$. By Lemma \ref{L2}, it follows that $H(u_1,:)\cdot{H(u_3,:)}=0$.
Then $\overline{h}_{11}h_{31}+\overline{h}_{12}h_{32}+\overline{h}_{13}h_{33}+\overline{h}_{14}h_{34}=-i-ih_{34}=0$.
Hence $h_{34}=-1$, which is a contradiction.

Subcase 2.2. $u_1\leftarrow{u_3}$ and $u_1\leftarrow{u_4}$.

By Lemma \ref{L2}, $H(u_1,:)\cdot{H(u_2,:)}=0$. Then $\overline{h}_{11}h_{21}+\overline{h}_{12}h_{22}+\overline{h}_{13}h_{23}+\overline{h}_{14}h_{24}=ih_{23}+ih_{24}=0$.
Hence $h_{23}=i, h_{24}=-i$ or $h_{23}=-i, h_{24}=i$, that is, $u_{2}\rightarrow{u_{3}},u_{2}\leftarrow{u_{4}}$ or $u_{2}\leftarrow{u_{3}},u_{2}\rightarrow{u_{4}}$. Without loss of generality, assume that $u_{2}\rightarrow{u_{3}},u_{2}\leftarrow{u_{4}}$. By Lemma \ref{L2}, it follows that $H(u_1,:)\cdot{H(u_3,:)}=0$.
Then $\overline{h}_{11}h_{31}+\overline{h}_{12}h_{32}+\overline{h}_{13}h_{33}+\overline{h}_{14}h_{34}=-i+ih_{34}=0$.
Hence $h_{34}=1$ i.e. there is an edge $u_{3}\leftrightarrow{u_{4}}$ in $G^{\phi}$. However, $H(u_1,:)\cdot{H(u_4,:)}=\overline{h}_{11}h_{41}+\overline{h}_{12}h_{42}+\overline{h}_{13}h_{43}+\overline{h}_{14}h_{44}=i+i\neq{0}$,  which is a contradiction.

Subcase 2.3. $u_1\rightarrow{u_3},\ u_1\leftarrow{u_4}$ or $u_1\leftarrow{u_3},\ u_1\rightarrow{u_4}$.

Without loss of generality, assume that $u_1\rightarrow{u_3}$ and $u_1\leftarrow{u_4}$.
By a similar way, we can prove that this subcase could not happen.

\textbf{Case 3.} No vertex has three incident edges, and there exists a vertex, say $u_1$,
has two incident edges  $u_1\leftrightarrow{u_2}$ and $u_1\leftrightarrow{u_3}$. Then for
the vertex $u_{4}$, there is an arc $u_{1}\rightarrow{u_{4}}$ or $u_{1}\leftarrow{u_{4}}$.

Suppose that $u_{1}\rightarrow{u_{4}}$. By Lemma \ref{L2}, we have $H(u_1,:)\cdot{H(u_2,:)}=0$.
Then $\overline{h}_{11}h_{21}+\overline{h}_{12}h_{22}+
\overline{h}_{13}h_{23}+\overline{h}_{14}h_{24}=h_{23}-ih_{24}=0$.
Hence $h_{23}=1, h_{24}=-i$ or $h_{23}=i, h_{24}=1$, that is,
$u_{2}\leftrightarrow{u_{3}},u_{2}\leftarrow{u_{4}}$ or $u_{2}\rightarrow{u_{3}},u_{2}\leftrightarrow{u_{4}}$.
If $u_{2}\leftrightarrow{u_{3}},u_{2}\leftarrow{u_{4}}$, then $H(u_1,:)\cdot{H(u_3,:)}=0$ from Lemma \ref{L2}.
It implies that $\overline{h}_{11}h_{31}+\overline{h}_{12}h_{32}+\overline{h}_{13}h_{33}+
\overline{h}_{14}h_{34}=1-ih_{34}=0$. Thus $h_{34}=-i$ i.e. $u_{3}\leftarrow{u_{4}}$. However, $H(u_1,:)\cdot{H(u_4,:)}=\overline{h}_{11}h_{41}+\overline{h}_{12}h_{42}+\overline{h}_{13}h_{43}+\overline{h}_{14}h_{44}=i+i\neq{0}$ , which is a contradiction. If $u_{2}\rightarrow{u_{3}},u_{2}\leftrightarrow{u_{4}}$, then $H(u_1,:)\cdot{H(u_3,:)}=0$
from Lemma \ref{L2}. It implies that $\overline{h}_{11}h_{31}+\overline{h}_{12}h_{32}+\overline{h}_{13}h_{33}+\overline{h}_{14}h_{34}=-i-ih_{34}=0$.
Thus $h_{34}=-1$, which is a contradiction.

For $u_{1}\leftarrow{u_{4}}$, we can prove that this case could not happen by the similar method.

\textbf{Case 4.} There exists a vertex, say $u_1$, has three incident edges
$u_1\leftrightarrow{u_2}$, $u_1\leftrightarrow{u_3}$ and $u_1\leftrightarrow{u_4}$.

Since $H(u_1,:)\cdot{H(u_2,:)}=0$, we can obtain that $h_{23}=i, h_{24}=-i$ or $h_{23}=-i, h_{24}=i$,
that is $u_{2}\rightarrow{u_{3}},u_{2}\leftarrow{u_{4}}$ or $u_{2}\leftarrow{u_{3}},u_{2}\rightarrow{u_{4}}$.
 Without loss of generality, assume that $u_{2}\rightarrow{u_{3}},u_{2}\leftarrow{u_{4}}$. Similarly,
 we have $h_{34}=i$ i.e. $u_{3}\rightarrow{u_{4}}$ by $H(u_1,:)\cdot{H(u_3,:)}=0$. That is
 $u_{2}\rightarrow{u_{3}}, u_{2}\leftarrow{u_{4}}$ and $u_{3}\rightarrow{u_{4}}$; see $G_1$ in Fig.\ref{Fig.3.3.}.

Thus, the proof is complete. \qed

Next, we determine all optimum Hermitian energy mixed graphs with underlying graph $Q_3$.

\begin{figure}[h,t,b,p]
\begin{center}
\scalebox{0.8}[0.8]{\includegraphics{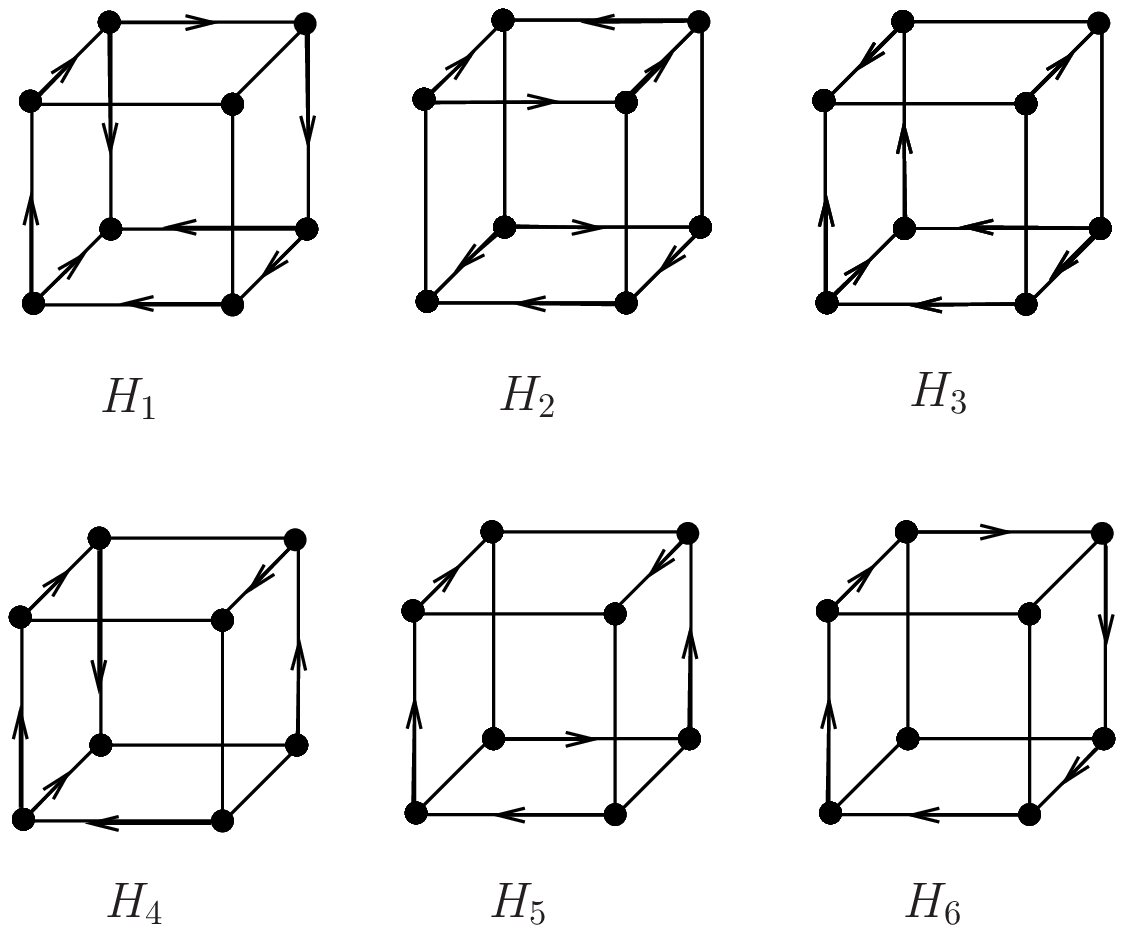}}
\end{center}
\caption{Optimum Hermitian energy mixed graphs with underlying graph $Q_3$.}\label{Fig.3.4.}
\end{figure}

\begin{thm}\label{T6}
Let $G^{\phi}$ be a $3$-regular optimum Hermitian energy mixed graph. If the underlying graph $G$ is $Q_3$, then $G^{\phi}$ (up to switching equivalence) is one of the following graphs: $D_2$ or $H_i$, where $i=1,2,...,6$; see Figs. \ref{Fig.2.2.} and \ref{Fig.3.4.}.
\end{thm}

\pf We divide our discussion into two cases:

\textbf{Case 1.} $G^{\phi}$ is an oriented graph.

From Lemma \ref{L5}, we obtain
that $G^{\phi}$ is $D_{2}$ drawn in Fig.\ref{Fig.2.2.}.

\begin{figure}[h,t,b,p]
\begin{center}
\scalebox{0.8}[0.8]{\includegraphics{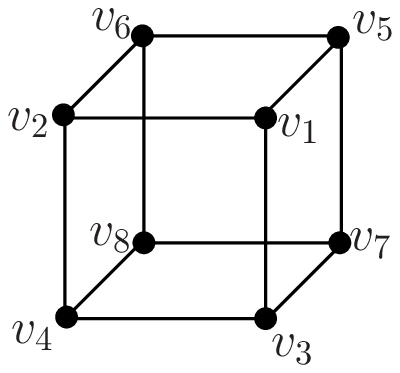}}
\end{center}
\caption{$Q_3$.}\label{Fig.3.5.}
\end{figure}

\textbf{Case 2.} $G^{\phi}$ is not an oriented graph. In the
following, we replace $G^{\phi}$ with ${Q_3^\phi}$ for convenience
and assume that $V(Q_{3})=\{v_1,v_2,...,v_8\}$, see
Fig.\ref{Fig.3.5.}.

Let $a$ (resp. $b$) denote the number of arcs (resp. undirected
edges) in $Q_{3}^{\phi}$, where $a+b=12$. Since $Q_{3}^{\phi}$ is
not an oriented graph, we get that $a\leq{11}$. Furthermore, there
are exactly six mixed $4$-cycles in $Q_{3}^{\phi}$. Let
$C_1,C_2,C_3,C_4,C_5$ and $C_6$ denote the mixed $4$-cycle induced
by vertices $\{v_1,v_2,v_4,v_3\}$, $\{v_2,v_4,v_8,v_6\}$,
$\{v_5,v_6,v_8,v_7\}$, $\{v_1,v_3,v_7,v_5\}$, $\{v_1,v_2,v_6,v_5\}$
and $\{v_3,v_4,v_8,v_7\}$, respectively. By Lemma \ref{L3}, we
deduce that every mixed $4$-cycle in $Q_{3}^{\phi}$ is one of the
three types in Fig.2.1. Thus, we obtain the following claim.

\textbf{Claim 1:}  In $Q_{3}^{\phi}$, every mixed $4$-cycle has
either two arcs and two undirected edges or four arcs.

It follows that each mixed $4$-cycle in $Q_{3}^{\phi}$ has at least
two arcs. Then, we have $a\geq\frac{2\times6}{2}=6$. Moreover, we
check that $a\neq{11,10}$. Consequently, $6\leq{a}\leq{9}$. Now we
divide the discussion about the values of $a$ and $b$ into four
subcases:

Subcase 2.1. $a=9$, $b=3$.

In this subcase, we want to determine three undirected edges in $Q_{3}^{\phi}$. Without loss of generality, suppose that $v_1\leftrightarrow{v_3}$. By Claim 1, both mixed $4$-cycle $C_1$ and $C_4$ have two undirected edges and hence we get the following four cases (up to isomorphism) by considering the other two undirected edges in $C_1$ and $C_4$.

(1) The other two undirected edges are $v_2\leftrightarrow{v_4}$ in $C_1$ and $v_5\leftrightarrow{v_7}$ in $C_4$. Then, there are three arcs in mixed $4$-cycle $C_2$ and $C_3$, which contradicts Claim 1.

(2) The other two undirected edges are $v_1\leftrightarrow{v_2}$ in $C_1$ and $v_5\leftrightarrow{v_7}$ in $C_4$. Then, there are three arcs in mixed $4$-cycle $C_5$ and $C_3$, which contradicts Claim 1.

(3) The other two undirected edges are $v_1\leftrightarrow{v_2}$ in $C_1$ and $v_3\leftrightarrow{v_7}$ in $C_4$. Then, there are three arcs in mixed $4$-cycle $C_5$ and $C_6$, which contradicts Claim 1.

(4) The other two undirected edges are $v_1\leftrightarrow{v_2}$ in $C_1$ and $v_1\leftrightarrow{v_5}$ in $C_4$. Then, mixed $4$-cycle $C_1$, $C_4$ and $C_5$ should be the first type in Fig.\ref{Fig.2.1.}; mixed $4$-cycle $C_2$, $C_3$ and $C_6$ should be the third type in Fig.\ref{Fig.2.1.}. Hence, there are two arcs $v_3\rightarrow{v_4},v_4\rightarrow{v_2}$ or $v_2\rightarrow{v_4},v_4\rightarrow{v_3}$ in $C_1$, $v_5\rightarrow{v_7},v_7\rightarrow{v_3}$ or $v_3\rightarrow{v_7},v_7\rightarrow{v_5}$ in $C_4$, and $v_2\rightarrow{v_6},v_6\rightarrow{v_5}$ or $v_5\rightarrow{v_6},v_6\rightarrow{v_2}$ in $C_5$. If there are two arcs $v_3\rightarrow{v_4}$ and $v_4\rightarrow{v_2}$ in a mixed graph, then we reverse every arc which is incident to vertex $v_4$ and acquire a new mixed graph. We can prove that the two mixed graphs are switching equivalent by the definition of switching equivalence. Without loss of generality, assume that $v_3\rightarrow{v_4},v_4\rightarrow{v_2}$. By a similar discussion, we assume that $v_2\rightarrow{v_6},v_6\rightarrow{v_5}$ and $v_5\rightarrow{v_7},v_7\rightarrow{v_3}$. Afterwards, we have either $v_4\rightarrow{v_8},v_6\rightarrow{v_8}$ or $v_8\rightarrow{v_6},v_8\rightarrow{v_4}$ in $C_2$. Analogously by switching equivalence, we assume that $v_4\rightarrow{v_8},v_6\rightarrow{v_8}$ and then $v_7\rightarrow{v_8}$. Therefore, we get the graph (up to switching equivalence) $H_1$ in Fig.\ref{Fig.3.4.}.

Subcase 2.2. $a=8$, $b=4$.

Now we want to determine four undirected edges in $Q_{3}^{\phi}$. Based on the discussion of subcase 2.1, we just need to find one more undirected edge.

If the three undirected edges which we have determined are $v_1\leftrightarrow{v_3},v_2\leftrightarrow{v_4}$ and $v_5\leftrightarrow{v_7}$, then there are two undirected edges in $C_1$ and $C_4$. Thus, the fourth undirected edge cannot be $v_1\leftrightarrow{v_2},v_3\leftrightarrow{v_4},v_1\leftrightarrow{v_5}$ or $v_3\leftrightarrow{v_7}$. If the fourth undirected edge is $v_2\leftrightarrow{v_6},v_5\leftrightarrow{v_6},v_4\leftrightarrow{v_8}$ or $v_7\leftrightarrow{v_8}$, then the resulting graphs are isomorphic. Without loss of generality, suppose that $v_2\leftrightarrow{v_6}$. Nevertheless, there are three arcs in $C_5$, which contradicts Claim 1. If the fourth undirected edge is $v_6\leftrightarrow{v_8}$, then mixed $4$-cycle $C_5$ and $C_6$ should be the third type in Fig.\ref{Fig.2.1.} and the others should be the second type in Fig.\ref{Fig.2.1.}. Thus, we get $H_2$ (up to isomorphism) in Fig.\ref{Fig.3.4.}.

If the three undirected edges which we have determined are $v_1\leftrightarrow{v_3},v_1\leftrightarrow{v_2}$ and $v_5\leftrightarrow{v_7}$, then there are two undirected edges in $C_1$ and $C_4$. Thus, the fourth undirected edge cannot be $v_2\leftrightarrow{v_4},v_3\leftrightarrow{v_4},v_1\leftrightarrow{v_5}$ or $v_3\leftrightarrow{v_7}$. If the fourth undirected edge is $v_2\leftrightarrow{v_6}$, then there are three arcs in $C_2$, which contradicts Claim 1. By a similar way, we deduce that the fourth undirected edge cannot be $v_6\leftrightarrow{v_8},v_4\leftrightarrow{v_8}$ or $v_7\leftrightarrow{v_8}$. If the fourth undirected edge is $v_5\leftrightarrow{v_6}$, then mixed $4$-cycle $C_1$ and $C_3$ should be the first type in Fig.\ref{Fig.2.1.}; mixed $4$-cycle $C_4$ and $C_5$ should be the second type in Fig.\ref{Fig.2.1.}; mixed $4$-cycle $C_2$ and $C_6$ should be the third type in Fig.\ref{Fig.2.1.}. Hence, there are two arcs $v_3\rightarrow{v_4},v_4\rightarrow{v_2}$ or $v_2\rightarrow{v_4},v_4\rightarrow{v_3}$ in $C_1$, and $v_7\rightarrow{v_8},v_8\rightarrow{v_6}$ or $v_6\rightarrow{v_8},v_8\rightarrow{v_7}$ in $C_3$. If there are two arcs $v_3\rightarrow{v_4}$ and $v_4\rightarrow{v_2}$ in a mixed graph, then we reverse every arc which is incident to vertex $v_4$ and acquire a new mixed graph. We can prove that the two mixed graphs are switching equivalent by the definition of switching equivalence. Without loss of generality, assume that $v_3\rightarrow{v_4},v_4\rightarrow{v_2}$. By a similar discussion, we assume that $v_7\rightarrow{v_8},v_8\rightarrow{v_6}$. Afterwards, we have either $v_4\rightarrow{v_8}$ or $v_8\rightarrow{v_4}$. If there is an arc $v_4\rightarrow{v_8}$, then we get the other arcs $v_6\rightarrow{v_2},v_1\rightarrow{v_5},v_7\rightarrow{v_3}$. Thus, we obtain $H_3$ depicted in Fig.\ref{Fig.3.4.}. If there is an arc $v_8\rightarrow{v_4}$, then we get the other arcs $v_2\rightarrow{v_6}$, $v_5\rightarrow{v_1}$, $v_3\rightarrow{v_7}$ and the resulting mixed graph is isomorphic to $H_3$.

If the three undirected edges which we have determined are $v_1\leftrightarrow{v_3},v_1\leftrightarrow{v_2}$ and $v_3\leftrightarrow{v_7}$, then there are two undirected edges in $C_1$ and $C_4$. Thus, the fourth undirected edge cannot be $v_2\leftrightarrow{v_4},v_3\leftrightarrow{v_4},v_1\leftrightarrow{v_5}$ or $v_5\leftrightarrow{v_7}$. If the fourth undirected edge is $v_2\leftrightarrow{v_6}$, then there are three arcs in $C_2$, which contradicts Claim 1. By a similar way, we deduce that the fourth undirected edge cannot be $v_5\leftrightarrow{v_6},v_6\leftrightarrow{v_8},v_4\leftrightarrow{v_8}$ or $v_7\leftrightarrow{v_8}$. Thus, this case could not happen.

If the three undirected edges which we have determined are $v_1\leftrightarrow{v_3},v_1\leftrightarrow{v_2}$ and $v_1\leftrightarrow{v_5}$, then there are two undirected edges in $C_1$, $C_4$ and $C_5$. Thus, the fourth undirected edge cannot be $v_2\leftrightarrow{v_4},v_3\leftrightarrow{v_4},v_3\leftrightarrow{v_7},v_5\leftrightarrow{v_7},v_2\leftrightarrow{v_6}$ or $v_5\leftrightarrow{v_6}$. If the fourth undirected edge is $v_4\leftrightarrow{v_8}$, then there are three arcs in $C_2$, which contradicts Claim 1. By a similar way, we deduce that the fourth undirected edge cannot be $v_6\leftrightarrow{v_8}$ or $v_7\leftrightarrow{v_8}$. Thus, this case could not happen.

Subcase 2.3. $a=7$, $b=5$.

Similarly in order to determine five undirected edges in $Q_{3}^{\phi}$, we just need to find two more undirected edges based on the discussion of subcase 2.1.

If the three undirected edges which we have determined are $v_1\leftrightarrow{v_3},v_2\leftrightarrow{v_4}$ and $v_5\leftrightarrow{v_7}$, then there are two undirected edges in $C_1$ and $C_4$. Thus, the other two undirected edges cannot be $v_1\leftrightarrow{v_2},v_3\leftrightarrow{v_4},v_1\leftrightarrow{v_5}$ or $v_3\leftrightarrow{v_7}$. If one of the other two undirected edges is $v_6\leftrightarrow{v_8}$, then there are two undirected edges in $C_2$ and $C_3$. By Claim 1, the last undirected edge cannot be $v_2\leftrightarrow{v_6},v_4\leftrightarrow{v_8},v_5\leftrightarrow{v_6}$ or $v_7\leftrightarrow{v_8}$. Then, there do not exist five undirected edges and hence this case could not happen. If there is an arc between $v_6$ and $v_8$, then the other two undirected edges (up to isomorphism) can be $v_2\leftrightarrow{v_6}$ and $v_4\leftrightarrow{v_8}$, $v_2\leftrightarrow{v_6}$ and $v_5\leftrightarrow{v_6}$, or $v_2\leftrightarrow{v_6}$ and $v_7\leftrightarrow{v_8}$. If the other two undirected edges are $v_2\leftrightarrow{v_6}$ and $v_4\leftrightarrow{v_8}$, then there are three undirected edges in $C_2$, which contradicts Claim 1. By a similar way, we deduce that the other two undirected edges cannot be $v_2\leftrightarrow{v_6}$ and $v_7\leftrightarrow{v_8}$. Therefore, the five undirected edges in $Q_{3}^{\phi}$ can be $v_1\leftrightarrow{v_3},v_2\leftrightarrow{v_4},v_5\leftrightarrow{v_7},v_2\leftrightarrow{v_6}$ and $v_5\leftrightarrow{v_6}$.

By a similar discussion, if the three undirected edges which we have determined are $v_1\leftrightarrow{v_3},v_1\leftrightarrow{v_2}$ and $v_5\leftrightarrow{v_7}$, then we deduce that the five undirected edges in $Q_{3}^{\phi}$ can be $v_1\leftrightarrow{v_3},v_1\leftrightarrow{v_2},v_5\leftrightarrow{v_7},v_2\leftrightarrow{v_6}$ and $v_6\leftrightarrow{v_8}$; if the three undirected edges which we have determined are $v_1\leftrightarrow{v_3},v_1\leftrightarrow{v_2}$ and $v_3\leftrightarrow{v_7}$, then we deduce that the five undirected edges in $Q_{3}^{\phi}$ can be either $v_1\leftrightarrow{v_3},v_1\leftrightarrow{v_2},v_3\leftrightarrow{v_7},v_2\leftrightarrow{v_6},v_4\leftrightarrow{v_8}$ or $v_1\leftrightarrow{v_3},v_1\leftrightarrow{v_2},v_3\leftrightarrow{v_7},v_5\leftrightarrow{v_6},v_7\leftrightarrow{v_8}$; if the three undirected edges which we have determined are $v_1\leftrightarrow{v_3},v_1\leftrightarrow{v_2}$ and $v_1\leftrightarrow{v_5}$, then this case could not happen. Regardless of the labels of vertices, the cases of the five undirected edges which we have determined are the same.

Without loss of generality, suppose that the five undirected edges are $v_1\leftrightarrow{v_3},v_1\leftrightarrow{v_2},v_3\leftrightarrow{v_7},v_5\leftrightarrow{v_6}$ and $v_7\leftrightarrow{v_8}$. Then, mixed $4$-cycle $C_1$, $C_4$ and $C_6$ should be the first type in Fig.\ref{Fig.2.1.}; mixed $4$-cycle $C_5$ and $C_3$ should be the second type in Fig.\ref{Fig.2.1.}; mixed $4$-cycle $C_2$ should be the third type in Fig.\ref{Fig.2.1.}. Hence, there are two arcs either $v_3\rightarrow{v_4},v_4\rightarrow{v_2}$ or $v_2\rightarrow{v_4},v_4\rightarrow{v_3}$ in $C_1$. If there are two arcs $v_3\rightarrow{v_4}$ and $v_4\rightarrow{v_2}$ in $C_1$, then we have an arc $v_4\rightarrow{v_8}$ in $C_6$. Otherwise, there is an arc $v_8\rightarrow{v_4}$. If there are three arcs $v_3\rightarrow{v_4},v_4\rightarrow{v_2}$ and $v_4\rightarrow{v_8}$ in a mixed graph, then we reverse every arc which is incident to vertex $v_4$ and acquire a new mixed graph. We can prove that the two mixed graphs are switching equivalent by the definition of switching equivalence. Without loss of generality, assume that $v_3\rightarrow{v_4},v_4\rightarrow{v_2}$ and $v_4\rightarrow{v_8}$. Afterwards, we have arcs either $v_7\rightarrow{v_5},v_5\rightarrow{v_1}$ or $v_1\rightarrow{v_5},v_5\rightarrow{v_7}$ in $C_4$. If there are two arcs $v_7\rightarrow{v_5}$ and $v_5\rightarrow{v_1}$ in $C_4$, then the other arcs are $v_2\rightarrow{v_6}$ and $v_6\rightarrow{v_8}$. Thus, we obtain $H_4$ depicted in Fig.\ref{Fig.3.4.}. If there are two arcs $v_1\rightarrow{v_5}$ and $v_5\rightarrow{v_7}$ in $C_4$, then the other arcs are $v_8\rightarrow{v_6},v_6\rightarrow{v_2}$ and the resulting mixed graph is isomorphic to $H_4$.

Subcase 2.4. $a=6$, $b=6$.

In order to determine six undirected edges in $Q_{3}^{\phi}$, we just need to find three more undirected edges based on the discussion of subcase 2.1.

If the three undirected edges which we have determined are $v_1\leftrightarrow{v_3},v_2\leftrightarrow{v_4}$ and $v_5\leftrightarrow{v_7}$, then there are two undirected edges in $C_1$ and $C_4$. Thus, the other three undirected edges cannot be $v_1\leftrightarrow{v_2},v_3\leftrightarrow{v_4},v_1\leftrightarrow{v_5}$ or $v_3\leftrightarrow{v_7}$. If one of the other three undirected edges is $v_2\leftrightarrow{v_6}$, then there must have an undirected edge $v_5\leftrightarrow{v_6}$ in $C_5$ by Claim 1. However, the last undirected edge cannot be $v_6\leftrightarrow{v_8},v_4\leftrightarrow{v_8}$ or $v_7\leftrightarrow{v_8}$ by Claim 1. Then, there do not exist six undirected edges. Similarly, we can show that one of the other three undirected edges cannot be $v_5\leftrightarrow{v_6},v_4\leftrightarrow{v_8}$ or $v_7\leftrightarrow{v_8}$. Hence, this case could not happen.

By a similar discussion, if the three undirected edges which we have determined are $v_1\leftrightarrow{v_3},v_1\leftrightarrow{v_2}$ and $v_5\leftrightarrow{v_7}$, then we deduce that the six undirected edges in $Q_{3}^{\phi}$ can be $v_1\leftrightarrow{v_3},v_1\leftrightarrow{v_2},v_5\leftrightarrow{v_7},v_2\leftrightarrow{v_6},v_4\leftrightarrow{v_8}$ and $v_7\leftrightarrow{v_8}$; if the three undirected edges which we have determined are $v_1\leftrightarrow{v_3},v_1\leftrightarrow{v_2}$ and $v_3\leftrightarrow{v_7}$, then we deduce that the six undirected edges in $Q_{3}^{\phi}$ can be either $v_1\leftrightarrow{v_3},v_1\leftrightarrow{v_2},v_3\leftrightarrow{v_7},v_2\leftrightarrow{v_6},v_6\leftrightarrow{v_8},v_7\leftrightarrow{v_8}$ or $v_1\leftrightarrow{v_3},v_1\leftrightarrow{v_2},v_3\leftrightarrow{v_7},v_5\leftrightarrow{v_6},v_6\leftrightarrow{v_8},v_4\leftrightarrow{v_8}$; if the three undirected edges which we have determined are $v_1\leftrightarrow{v_3},v_1\leftrightarrow{v_2}$ and $v_1\leftrightarrow{v_5}$, then we deduce that the six undirected edges in $Q_{3}^{\phi}$ can be $v_1\leftrightarrow{v_3},v_1\leftrightarrow{v_2},v_1\leftrightarrow{v_5},v_6\leftrightarrow{v_8},v_4\leftrightarrow{v_8}$ and $v_7\leftrightarrow{v_8}$. Regardless of the labels of vertices, the case that $v_1\leftrightarrow{v_3},v_1\leftrightarrow{v_2},v_5\leftrightarrow{v_7},v_2\leftrightarrow{v_6},v_4\leftrightarrow{v_8},v_7\leftrightarrow{v_8}$ is the same with the case that $v_1\leftrightarrow{v_3},v_1\leftrightarrow{v_2},v_3\leftrightarrow{v_7},v_5\leftrightarrow{v_6},v_6\leftrightarrow{v_8},v_4\leftrightarrow{v_8}$.
Thus, we get the following three cases.

(1) The six undirected edges are $v_1\leftrightarrow{v_3},v_1\leftrightarrow{v_2},v_3\leftrightarrow{v_7},v_2\leftrightarrow{v_6},v_6\leftrightarrow{v_8}$ and $v_7\leftrightarrow{v_8}$. Then, every mixed $4$-cycle in $Q_{3}^{\phi}$ should be the first type in Fig.\ref{Fig.2.1.}. Hence, there are two arcs either $v_3\rightarrow{v_4},v_4\rightarrow{v_2}$ or $v_2\rightarrow{v_4},v_4\rightarrow{v_3}$ in $C_1$. If there are two arcs $v_3\rightarrow{v_4}$ and $v_4\rightarrow{v_2}$ in $C_1$, then we get an arc $v_4\rightarrow{v_8}$ in $C_6$ and $v_8\rightarrow{v_4}$ in $C_2$, a contradiction. Analogously, the case that there are two arcs $v_2\rightarrow{v_4}$ and $v_4\rightarrow{v_3}$ in $C_1$ could not happen.

(2) The six undirected edges are $v_1\leftrightarrow{v_3},v_1\leftrightarrow{v_2},v_3\leftrightarrow{v_7},v_5\leftrightarrow{v_6},v_6\leftrightarrow{v_8}$ and $v_4\leftrightarrow{v_8}$. Then, mixed $4$-cycle $C_1$, $C_2$, $C_3$ and $C_4$ should be the first type in Fig.\ref{Fig.2.1.}; mixed $4$-cycle $C_5$ and $C_6$ should be the second type in Fig.\ref{Fig.2.1.}. Hence, there are two arcs either $v_3\rightarrow{v_4},v_4\rightarrow{v_2}$ or $v_2\rightarrow{v_4},v_4\rightarrow{v_3}$ in $C_1$. If there are two arcs $v_3\rightarrow{v_4}$ and $v_4\rightarrow{v_2}$ in $C_1$, then we get the other arcs $v_2\rightarrow{v_6},v_5\rightarrow{v_1},v_7\rightarrow{v_5}$ and $v_8\rightarrow{v_7}$. Thus, we obtain $H_5$ depicted in Fig.\ref{Fig.3.4.}. If there are two arcs $v_2\rightarrow{v_4}$ and $v_4\rightarrow{v_3}$ in $C_1$, then we get the other arcs $v_6\rightarrow{v_2},v_1\rightarrow{v_5},v_5\rightarrow{v_7}$, $v_7\rightarrow{v_8}$ and the resulting mixed graph is isomorphic to $H_5$.

(3) The six undirected edges are $v_1\leftrightarrow{v_3},v_1\leftrightarrow{v_2},v_1\leftrightarrow{v_5},v_6\leftrightarrow{v_8},v_4\leftrightarrow{v_8}$ and $v_7\leftrightarrow{v_8}$. Then, every mixed $4$-cycle in $Q_{3}^{\phi}$ should be the first type in Fig.\ref{Fig.2.1.}. Hence, there are two arcs either $v_3\rightarrow{v_4},v_4\rightarrow{v_2}$ or $v_2\rightarrow{v_4},v_4\rightarrow{v_3}$ in $C_1$. If there are two arcs $v_3\rightarrow{v_4}$ and $v_4\rightarrow{v_2}$ in $C_1$, then we get the other arcs $v_2\rightarrow{v_6},v_6\rightarrow{v_5},v_5\rightarrow{v_7}$ and $v_7\rightarrow{v_3}$. Thus, we obtain $H_6$ depicted in Fig.\ref{Fig.3.4.}. If there are two arcs $v_2\rightarrow{v_4}$ and $v_4\rightarrow{v_3}$ in $C_1$, then we get the other arcs $v_6\rightarrow{v_2},v_5\rightarrow{v_6},v_7\rightarrow{v_5},v_3\rightarrow{v_7}$ and the resulting mixed graph is isomorphic to $H_6$.

Thus, the proof is complete. \qed

\section{The uniqueness of oriented graph $Q_k^{\phi}$ with optimum Hermitian energy}

It is difficult to determine all optimum Hermitian energy mixed graphs with underlying graph hypercube $Q_k$ for $k\geq{4}$. In \cite{T}, Tian gave $Q_k$ an orientation such that it has optimum Hermitian energy. Besides, Gong and Xu \cite{2} proved that $3$-regular optimum Hermitian energy oriented graph $Q_{3}^{\phi}$ is unique (up to switching equivalence). In this section, we show that any optimum Hermitian energy oriented graph with underlying graph $Q_k$ is unique (up to switching equivalence) for any positive integer $k$.

Firstly, we give the following definition about hypercube $Q_{k}$ which can be found in \cite{4}.

\begin{defi}\cite{4}\label{D6}
A hypercube $Q_{k}$ of dimension $k$ is defined recursively in terms of the Cartesian
product of graphs as follows
\begin{equation*}
Q_{k}=
\begin{cases}
K_{2}, & k=1,\\
Q_{k-1}\Box{Q_{1}}, & k\geq{2}.
\end{cases}
\end{equation*}
\end{defi}

\begin{lem}\label{L3.2}
Let $Q_{k}^{\phi}$ be an oriented graph with the orientation $\phi$.
Then $Q_{k}^{\phi}$ has optimum Hermitian energy if and only if
every mixed $4$-cycle in $Q_{k}^{\phi}$ is the third type in
Fig.\ref{Fig.2.1.}.
\end{lem}

\pf For any two distinct vertices $u$ and $v$ of $Q_{k}$, we know that $|N(u)\cap{N(v)}|$
is either zero or two, where $N(\cdot)$ stands for the neighborhood of a vertex in $Q_{k}$.
Thus, if there is one common neighbor between the two vertices in $Q_{k}^{\phi}$, then they
have exactly two common neighbors $x$ and $y$ i.e. there is exactly one mixed $4$-cycle $uxvy$.
By Lemma \ref{L3}, it is easy to obtain this lemma. \qed

Hypercube $Q_{k}$ is a very important family of graphs and it has
many nice properties. Liu and Li \cite{1} gave $Q_{k}$ a suitable
orientation such that it has optimum Hermitian energy. Now we give
$Q_{k}$ a new orientation $\phi{_{0}}$. For convenience, we assume
that the vertex set of $Q_{k}$ is
$\{1,2,\ldots,2^{k-1},2^{k-1}+1,2^{k-1}+2,\ldots,2^{k}\}$ with
$G[V_{1}]=G[V_{2}]=Q_{k-1}$, where $V_{1}=\{1,2,\ldots,2^{k-1}\}$,
$V_{2}=\{2^{k-1}+1,2^{k-1}+2,\ldots,2^{k}\}$. Firstly, we give the
hypercube $Q_{1}$ an orientation $Q_{1}^{\phi{_{0}}}$ such that
$1\rightarrow{2}$. Afterwards, we suppose that $Q_{k-1}$ has been
oriented into $Q_{k-1}^{\phi{_{0}}}$. By reversing every arc of
$Q_{k-1}^{\phi{_{0}}}$, we can get another new orientation denoted
by $-\phi{_{0}}$. For $Q_{k}$, we give $G[V_{1}]$ the orientation
${\phi{_{0}}}$ and $G[V_{2}]$ the orientation ${-\phi{_{0}}}$. Next
we put an arc from each vertex in $G[V_{1}]$ to the corresponding
vertex in $G[V_{2}]$ i.e. $t\rightarrow{2^{k-1}+t}$ for
$t=1,2,\ldots,2^{k-1}$. Then we get $Q_{k}^{\phi{_{0}}}$.

The lemma below shows that the Hermitian energy of $Q_{k}^{\phi{_{0}}}$ is optimum.

\begin{lem}\label{L8}
Let $Q_{k}$ be a hypercube of dimension $k$ with $n=2^{k}$ vertices.
Then $Q_{k}^{\phi{_{0}}}$ satisfies $H^{2}(Q_{k}^{\phi{_{0}}})=kI_{n}$
(or $\mathcal{E}_H(Q_{k}^{\phi{_{0}}})={n\sqrt{k}}$).
\end{lem}

\pf If $k=1$, then
\[
H(Q_{1}^{\phi{_{0}}})=\begin{bmatrix}
0 & i \\-i & 0\\
\end{bmatrix}.
\]
Hence it is easy to show that $H^{2}(Q_{1}^{\phi{_{0}}})=I_{2}$ and we only
need to consider the case $k\geq{2}$.

By Lemma \ref{L3.2}, if every mixed $4$-cycle of $Q_{k}^{\phi{_{0}}}$ is the
third type in Fig.\ref{Fig.2.1.}, then $\mathcal{E}_H(Q_{k}^{\phi{_{0}}})={n\sqrt{k}}$.
Therefore, we just need to show that every mixed $4$-cycle of $Q_{k}^{\phi{_{0}}}$ is
the third type in Fig.\ref{Fig.2.1.} for $k\geq{2}$. We shall apply induction on $k$.
If $k=2$, $Q_{2}^{\phi{_{0}}}$ is the third type in Fig.\ref{Fig.2.1.}. Suppose now that
$k>2$ and the lemma holds for fewer $k$. Then every mixed $4$-cycle of $Q_{k-1}^{\phi{_{0}}}$
is the third type in Fig.\ref{Fig.2.1.} and so is $Q_{k-1}^{-\phi{_{0}}}$. Moreover by the
definition of $Q_{k}^{\phi{_{0}}}$, we have
\[
H(Q_{k}^{\phi{_{0}}})=\begin{bmatrix}
H(Q_{k-1}^{\phi{_{0}}}) & iI \\-iI & H(Q_{k-1}^{-\phi{_{0}}})\\
\end{bmatrix}.
\]
For any mixed $4$-cycle of $Q_{k}^{\phi{_{0}}}$, if it is contained
by the induced subgraph $G[V_1]$ or $G[V_2]$, then we get that it is
the third type in Fig.\ref{Fig.2.1.} by the induction hypothesis.
Thus, we just talk about the mixed $4$-cycle $C$ induced by
$\{s,t,t+2^{k-1},s+2^{k-1}\}$, where $1\leq{s},t\leq{2^{k-1}}$ and
$s\neq{t}$. Without loss of generality, suppose that there is an arc
$s\rightarrow{t}$ in $Q_{k}^{\phi{_{0}}}$. Then, we have arcs
${s+2^{k-1}}\leftarrow{t+2^{k-1}},s\rightarrow{s+2^{k-1}}$ and
$t\rightarrow{t+2^{k-1}}$. It is easy to see that the mixed
$4$-cycle $C$ is the third type in Fig.\ref{Fig.2.1.}. Above all, we
get that every mixed $4$-cycle of $Q_{k}^{\phi{_{0}}}$ is the third
type in Fig.\ref{Fig.2.1.} for $k\geq{2}$. Thus, We complete the
proof. \qed

\begin{thm}\label{T8}
Let $Q_{k}^{\phi}$ be an optimum Hermitian energy mixed graph with
underlying graph $Q_{k}$. If $Q_{k}^{\phi}$ is an oriented graph,
then $Q_{k}^{\phi}$ is unique (up to switching equivalence) for any
positive integer $k$.
\end{thm}

\pf We shall apply induction on $k$. If $k=1$, $Q_{k}^{\phi}$ is an
arc. Thus, we obtain that $Q_{k}^{\phi}$ is unique (up to switching
equivalence) for $k=1$. Now we assume that the theorem holds for
fewer $k$.

Let $Q_{k}^{\phi}$ be an oriented graph with optimum Hermitian
energy. For the sake of convenience, assume that
$V(Q_{k}^{\phi})=\{1,2,\ldots,2^{k-1},2^{k-1}+1,2^{k-1}+2,\ldots,2^{k}\}$
such that both $G[V_{1}]$ and $G[V_{2}]$ are mixed graphs with
underlying graph $Q_{k-1}$, where $V_{1}=\{1,2,\ldots,2^{k-1}\}$ ,
$V_{2}=\{2^{k-1}+1,2^{k-1}+2,\ldots,2^{k}\}$. From Lemma \ref{L3.2},
we know that every mixed 4-cycle in $Q_{k}^{\phi}$ is the third type
in Fig.\ref{Fig.2.1.}. Then every mixed 4-cycle in $G[V_{1}]$ and
$G[V_{2}]$ is the third type in Fig.\ref{Fig.2.1.}. By Lemma
\ref{L3.2}, $G[V_{1}]$ and $G[V_{2}]$ have optimum Hermitian energy.
Thus, we have
\[
H(Q_{k}^{\phi})=\begin{bmatrix}
H(G[V_{1}]) & S \\S^{*} & H(G[V_{2}])\\
\end{bmatrix}\ ,
\]
where $S$ is a diagonal matrix and each diagonal element belongs to
$\{i,-i\}$. By the induction hypothesis, we can find a switching
matrix $D_{1}(\theta)$ such that
$H(Q_{k-1}^{\phi_{0}})=D_{1}^{-1}(\theta)H(G[V_{1}])D_{1}(\theta)$.
Assume that $Q_{k-1}$ is a bipartite graph with bipartition $X$ and
$Y$. Let a diagonal matrix
$D_{3}(\theta)=diag(\theta{(v_{k})}|\theta{(v_{k})}=1 ,\ if\
v_{k}\in{X};\ \theta{(v_{k})}=-1, \ if\ v_{k}\in{Y})$. Then
$H(Q_{k-1}^{-\phi_{0}})=D_{3}^{-1}(\theta)H(Q_{k-1}^{\phi_{0}})D_{3}(\theta)$.
Hence $Q_{k-1}^{-\phi_{0}}\sim{Q_{k-1}^{\phi_{0}}}$ i.e.
$Q_{k-1}^{-\phi_{0}}$ has optimum Hermitian energy. Similarly, we
can find the switching matrix $D_{2}(\theta)$ such that
$H(Q_{k-1}^{-\phi_{0}})=D_{2}^{-1}(\theta)H(G[V_{2}])D_{2}(\theta)$
by the induction hypothesis.

Let \[
T_{1}=\begin{bmatrix}
D_{1}(\theta) & 0 \\0 & I\\
\end{bmatrix}\ ,\ \ \
T_{2}=\begin{bmatrix}
I & 0 \\0 & D_{2}(\theta)\\
\end{bmatrix}.
\]
Then $T_{1}$ and $T_{2}$ are switching matrices and we get that

\begin{eqnarray*}
&&T_2^{-1}T_1^{-1}H(Q_{k}^{\phi})T_{1}T_{2}\\
&=&{\begin{bmatrix}
I & 0 \\0 & D_{2}(\theta)\\
\end{bmatrix}}^{-1}
\begin{bmatrix}
D_{1}(\theta) & 0 \\0 & I\\
\end{bmatrix}^{-1}
\begin{bmatrix}
H(G[V_{1}]) & S \\S^{*} & H(G[V_{2}])\\
\end{bmatrix}
\begin{bmatrix}
D_{1}(\theta) & 0 \\0 & I\\
\end{bmatrix}
\begin{bmatrix}
I & 0 \\0 & D_{2}(\theta)\\
\end{bmatrix}\\
&=&\begin{bmatrix}
D_{1}(\theta)^{-1}H(G[V_{1}])D_{1}(\theta) & D_{1}
(\theta)^{-1}SD_{2}(\theta) \\D_{2}(\theta)^{-1}S^{*}D_{1}
(\theta) & D_{2}(\theta)^{-1}H(G[V_{2}])D_{2}(\theta)\\
\end{bmatrix}\\
&=&\begin{bmatrix}
H(Q_{k-1}^{\phi_{0}}) & S_{1} \\S_{1}^{*} & H(Q_{k-1}^{-\phi_{0}})\\
\end{bmatrix},
\end{eqnarray*}

where $S_{1}=D_{1}^{-1}(\theta)SD_{2}(\theta)=\begin{bmatrix}
s_{11} & 0 & \cdots & 0 \\0 & s_{22} & \cdots &0 \\
\vdots & \vdots & \ddots & \vdots \\0 & 0 & \cdots & s_{2^{k-1}2^{k-1}} \\
\end{bmatrix}\ $ and $s_{tt}\in{\\\\\{i,-i\}}$ with $1\leq{t}\leq{2^{k-1}}$.

Now we divide the discussion about the value of $s_{11}$ into two cases:

\textbf{Case 1.} $s_{11}=i$.

Let $T_{3}=I$. Then
\[
H=T_{3}^{-1}T_{2}^{-1}T^{-1}_{1}H(Q_{k}^{\phi})T_{1}T_{2}T_{3}=\begin{bmatrix}
H(Q_{k-1}^{\phi_{0}}) & S_{2} \\S_{2}^{*} & H(Q_{k-1}^{-\phi_{0}})\\
\end{bmatrix},
\]
where $S_{2}=I^{-1}S_{1}(I)=S_{1}$.

\textbf{Case 2.} $s_{11}=-i$.

Let
\[
T_{3}=\begin{bmatrix}
I & 0 \\0 & -I\\
\end{bmatrix}.
\] Then
\[
H=T_{3}^{-1}T_{2}^{-1}T^{-1}_{1}H(Q_{k}^{\phi})T_{1}T_{2}T_{3}=\begin{bmatrix}
H(Q_{k-1}^{\phi_{0}}) & S_{2} \\S_{2}^{*} & H(Q_{k-1}^{-\phi_{0}})\\
\end{bmatrix},
\]
where $S_{2}=I^{-1}S_{1}(-I)=-S_{1}$.

\begin{figure}[h,t,b,p]
\begin{center}
\scalebox{0.8}[0.8]{\includegraphics{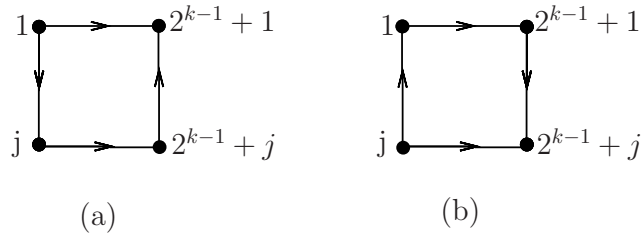}}
\end{center}
\caption{Two orientations of edges related to vertices $1,j,2^{k-1}+1,2^{k-1}+j$.}\label{Fig.4.6.}
\end{figure}

Assume that $S_{2}=\begin{bmatrix}
s_{11}^{'} & 0 & \cdots & 0 \\0 & s_{22}^{'} & \cdots &0 \\
\vdots & \vdots & \ddots & \vdots \\0 & 0 & \cdots & s_{2^{k-1}2^{k-1}}^{'} \\
\end{bmatrix}$ and $H=(h_{kl})$. Then $h_{t(2^{k-1}+t)}=s^{'}_{tt},1\leq{t}\leq{2^{k-1}}$
and $h_{1(2^{k-1}+1)}=s_{11}^{'}=i$. Let $T=T_{1}T_{2}T_{3}$. Then
$H=T^{-1}H(Q_{k}^{\phi})T$, where $T$ is a diagonal matrix and every
diagonal element belongs to $\{1,-1\}$. Since
$H^{2}(Q_{k}^{\phi})=kI$, we get that $HH^*=H^{2}=kI$. Thus, the
inner product of any two rows of $H$ is zero. Suppose that vertex
$j$ is a neighbor of vertex $1$, where $j\in{V_{1}}$. Next we
consider the inner product of the first row and the $2^{k-1}+j$th
row in $H$. If $1\rightarrow{j}$, then
$h_{1j}=i,h_{(2^{k-1}+1)(2^{k-1}+j)}=-i$ and
$h_{(2^{k-1}+j)(2^{k-1}+1)}=\overline{h}_{(2^{k-1}+1)(2^{k-1}+j)}=i$.
Since
$H(1,:)\cdot{H(2^{k-1}+j,:)}=h_{1j}\overline{h}_{(2^{k-1}+j)j}+h_{1(2^{k-1}+1)}
\overline{h}_{(2^{k-1}+j)(2^{k-1}+1)}=ih_{j(2^{k-1}+j)}+i(-i)=0$, we
get that $h_{j(2^{k-1}+j)}=i$, i.e., $j\rightarrow{2^{k-1}+j}$; see
Fig.\ref{Fig.4.6.}(a). If $1\leftarrow{j}$, then
$h_{1j}=-i,h_{(2^{k-1}+1)(2^{k-1}+j)}=i$ and
$h_{(2^{k-1}+j)(2^{k-1}+1)}=\overline{h}_{(2^{k-1}+1)(2^{k-1}+j)}=-i$.
Since
$H(1,:)\cdot{H(2^{k-1}+j,:)}=h_{1j}\overline{h}_{(2^{k-1}+j)j}+h_{1(2^{k-1}+1)}
\overline{h}_{(2^{k-1}+j)(2^{k-1}+1)}=(-i)h_{j(2^{k-1}+j)}+ii=0$, we
get that $h_{j(2^{k-1}+j)}=i$, i.e., $j\rightarrow{2^{k-1}+j}$; see
Fig.\ref{Fig.4.6.}(b). Due to the connection of $Q_{k-1}$, we can
show that
$h_{2(2^{k-1}+2)}=h_{3(2^{k-1}+3)}=\cdots{=}h_{2^{k-1}2^{k}}=i$.
Then $S_{2}=iI$ and $H=
\begin{bmatrix}
H(Q_{k-1}^{\phi_{0}}) & iI \\-iI & H(Q_{k-1}^{-\phi_{0}})\\
\end{bmatrix}=H(Q_{k}^{\phi_{0}})$. Thus, there exists a switching matrix
$D(\theta)=T$ such that $H(Q_{k}^{\phi_{0}})=D(\theta)^{-1}H(Q_{k}^{\phi})D(\theta)$.
That is $Q_{k}^{\phi}\sim{Q_{k}^{\phi_{0}}}$.

Above all, we conclude that any optimum Hermitian energy oriented graph with
underlying graph $Q_{k}$ is unique (up to switching equivalence)
for any positive integer $k$. The proof is complete.\qed

\noindent\textbf{Remark 4.1} The optimum Hermitian energy orientation ${\phi{_{0}}}$ of $Q_{k}$ is similar to the orientation obtained by Tian \cite{T}. However, our proof is very different from Tian's. His proof is based on the skew adjacency matrix, while ours uses the third type mixed 4-cycle. Moreover, we prove that the optimum Hermitian energy orientations of $Q_{k}$ are unique (up to switching equivalence).

\end{document}